\documentclass{article}
\usepackage{amssymb,amsmath,graphicx}
\newcommand{\R}{{\mathbb R}}
\newcommand{\eone}{e_1}
\newcommand{\sgn}{\mathrm{sgn}}
\newcommand{\matlab}{{\sc matlab}}
\newcommand{\epsmch}{\epsilon_{\mathrm{mch}}}
\title{On the Choice of Sign Defining\\ Householder Transformations}
\author{
Michael L. Overton\thanks{Courant Institute of Mathematical Sciences, New York University, {\tt mo1@nyu.edu}} 
\and Pinze Yu\thanks{Courant Institute of Mathematical Sciences, New York University, {\tt py2050@nyu.edu}} 
}
\date{October 7, 2023}
\begin{document}
\maketitle
\begin{abstract}
It is well known that, when defining Householder transformations, the correct choice of sign in the
standard formula is important to avoid cancellation and hence numerical instability. 
In this note we point out that when the ``wrong'' choice of sign is used,
the extent of the resulting instability depends in a somewhat subtle way 
on the data leading to cancellation.
\end{abstract}

\small
\centerline{AMS Subject Classification: 65F05}

\normalsize
\section{Introduction}
The QR factorization is a standard tool in numerical linear algebra, and Householder transformations
provide the best general method to compute it. Following \cite[Sec.~19.1]{Hig02}, a Householder transformation (or
Householder reflector) has the form
\begin{equation}\label{Pdef}
      P = I -\frac{2}{v^{T}v} vv^{T},
\end{equation}
where $I$ is the identity matrix and $v$ is a nonzero vector.
It is easily verified that $P$ is an orthogonal matrix, i.e., $P^{T}P=I$. 
The first step in the Householder reduction of an $m\times n$ matrix $A$, with $m\geq n$,
to triangular form is to define a Householder transformation $P_{1}$ that maps $x$, the first column of $A$,
to a multiple of the first coordinate vector $\eone=[1,0,\ldots,0]^{T} \in\R^{m}$. 
Since $P_{1}x$ must have the same Euclidean length
as $x$, we require $P_{1}x=\sigma \|x\|\eone$, where $\sigma=\pm 1$ and $\|\cdot\|$ denotes
the 2-norm. Thus we need
$$
    P_{1}x = x - \frac{2v^{T}x}{v^{T}v} v = \sigma \|x\|\eone
$$
which implies that $v$ is a scalar multiple of $x - \sigma \|x\| \eone$, and since $P_{1}$
is independent of $\|v\|$, without loss of generality we can choose
\begin{equation}\label{vdef}
       v = x - \sigma \|x\| \eone.
\end{equation}
To avoid numerical cancellation in \eqref{vdef}, it is generally recommended to use 
\begin{equation}\label{right}
      \sigma = - \sgn(x_{1})
\end{equation}
where $x_{1}$ is the first component of the vector $x$ and $\sgn$
is the standard sign function, which for convenience we define to be $+1$ if its argument is 
zero. The transformation $P_{1}$ is then applied to the remaining columns of $A$ as well, exploiting the
formula \eqref{Pdef} for efficiency, yielding the matrix $P_{1}A$ whose first column has 
all zeros except in the first position.
The factorization is completed by repeating the process for
every column of $A$, working with only with the data in rows
$k$ through $m$ and columns $k$ through $n$ at the $k$th step, yielding a total of~$n$ Householder transformations $P_{1},P_{2}\ldots,P_{n}$,
along with the upper triangular final matrix $R$. Then in exact arithmetic, $A=QR$,
with $Q=P_{1}P_{2}\ldots P_{n}$.

In this note we examine exactly what occurs if the ``wrong'' sign\footnote{It is pointed out in
\cite[Sec.~19.1]{Hig02} that the sign \eqref{wrong} may be used if
the formula for $v$ is rearranged; see \cite{Par71} for details.
While this is useful if 
consistent signs are preferred in computing the QR factorization, it is not relevant to the subsequent discussion.}
\begin{equation}\label{wrong}
      \sigma = \sgn(x_{1})
\end{equation}
is used to compute $v$ in \eqref{vdef}.

\section{Observation}

%For some years the first author has asked students in his undergraduate numerical computing class
%to construct a numerical experiment showing that using the wrong sign \eqref{wrong}
%in the Householder reduction
%gives an unstable result. 
We consider the following experiment. We would like to choose
$A$ so that using the wrong sign \eqref{wrong} results in as much cancellation as possible;
an easy way to do this is to choose the first column to have much smaller entries, 
in magnitude, than the (1,1) entry, so that $\sgn(x_{1}) \|x\|$ approximately cancels with $x_{1}$ in \eqref{vdef}. Here, we report the
results of an experiment computing $Q$ and $R$ using both choices of sign 
for a $3\times 2$ matrix $A$ with $a_{11} = 1$, $a_{21}=\delta$, $a_{32}=0$ and the second column chosen randomly,
for $\delta$ taking the successive
values $10^{-p}$, $p=1,2,\ldots,16$. 
The experiment was conducted using \matlab\ on a MacBook Pro, 
for which the machine epsilon $\epsmch$ (the gap between 1 and the next larger floating
point number) is approximately $10^{{-16}}$ (as \matlab\ uses IEEE double precision
by default).

\begin{figure}[h]
\begin{center}
\includegraphics[scale=0.55]{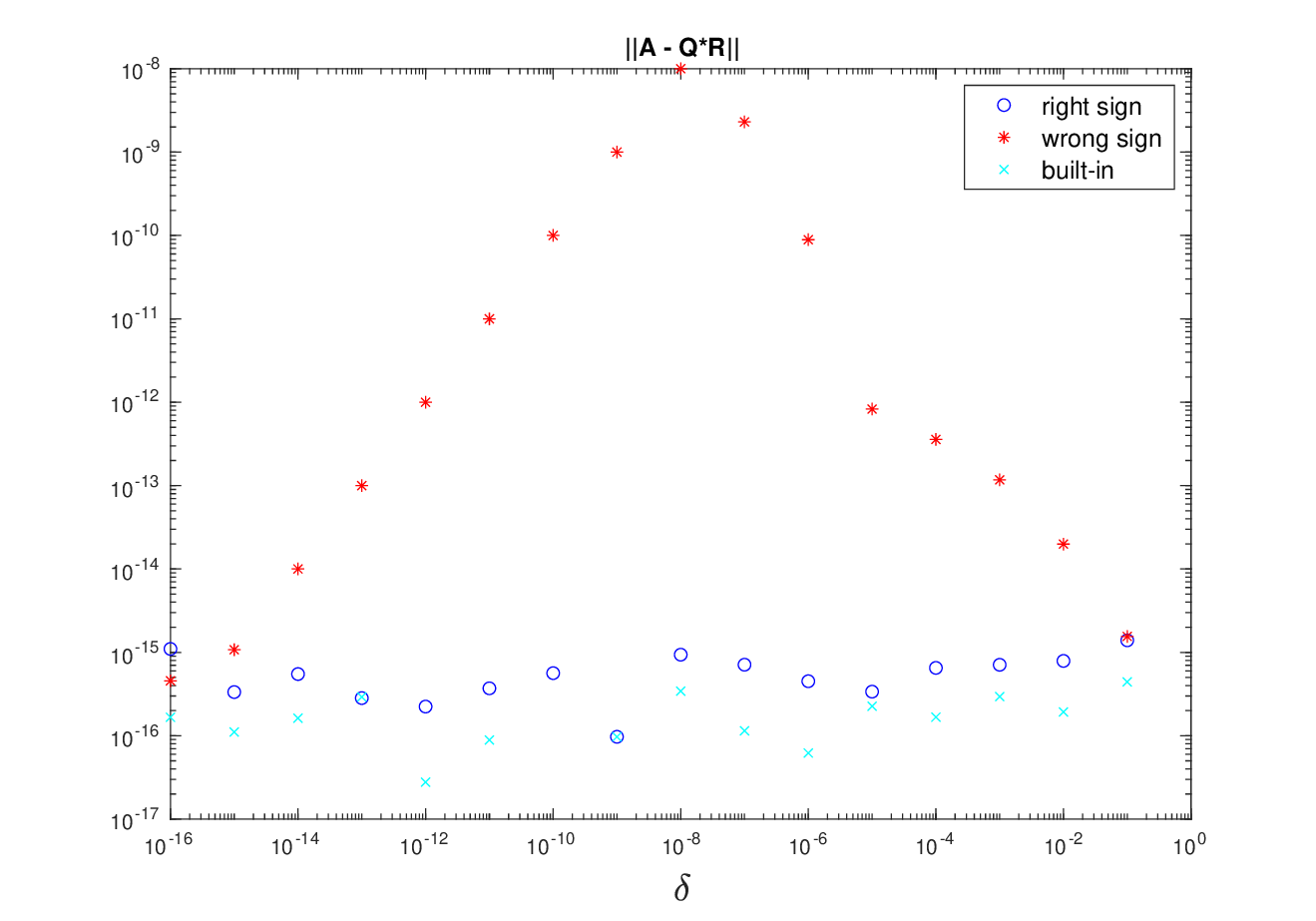}
\end{center}
\caption{The 2-norm of $A-QR$, where $Q$ and $R$ are the computed Q and~R factors of a
$3\times 2$ matrix $A$
with first column $[1,~\delta,~0]^{T}$, using Householder reduction with
the correct choice of sign \eqref{right} (blue circles), the wrong choice of sign \eqref{wrong}
(red asterisks) and \matlab's built-in {\tt qr} function (cyan crosses), all plotted
as a function of $\delta$.\label{fig1}}
\end{figure}

Figure \ref{fig1} shows the computed 2-norm
$\|A - QR\|$ for each choice of $\delta$ and for three algorithms: using the
correct sign (blue circles), the wrong sign (red asterisks), and using \matlab's built-in \verb@qr@ (cyan crosses); note the log-log scaling. 
Unsurprisingly, the results using the correct
choice of sign or the built-in \verb@qr@ are, for all $\delta$, approximately 
$\epsmch$. Surprisingly, however, the results using
the wrong sign appear in an inverted-V pattern with respect to $\delta$.
This is somewhat reminiscent of the well-known
V pattern that is often used, for example in \cite[Chap.~11]{Ove01}, to show how
the truncation error and rounding error respectively dominate the error in
the approximation of a derivative of a function~$f$ at a point $x$ 
by a finite difference quotient $\frac{f(x+h)-f(x)}{h}$, the former dominant
for large $h$ and the latter dominant for small $h$. The comparison
even extends to noting that the right side of the inverted V is ragged, indicating
dominance by rounding error, while the left side is a straight line, indicating purely
linear dependence; in the finite-difference example, the roles of left and right are reversed. 
Note that the choice of $\delta\approx \epsmch^{1/2}$, the square root of the machine precision, gives the most inaccurate result, while in the finite difference example,
 it is well known that $h\approx \epsmch^{1/2}$ is the best choice, assuming appropriately scaled data.
The results shown in Figure~\ref{fig1} are essentially unchanged if much larger matrices are
used.

\section{Explanation}

The right side of the inverted V, where
the error increases as $\delta$ decreases, is what we expected as the cancellation error in
\eqref{vdef} becomes more dominant. But what about the left side, where the error decreases
as $\delta$ continues to decrease? In fact, this is easily explained.
In the experiment, the first
column of $A$ is $[1,~\delta,~0]^{T}$, whose 2-norm is $\sqrt{1 + \delta^{2}}$, so for $\delta$
somewhat less than $\epsmch^{1/2}$, the computed 2-norm is precisely 1. This results
in the first component of the vector $v$ defining the first Householder transformation
being zero. The second component of $v$ is $\delta$ and the third is zero, 
so the normalized vector $v/\|v\|$
is the second unit vector. This means that the first Householder transformation is the identity
except with $-1$ instead of $+1$ in the (2,2) position. Thus the first column of $Q$,
the product of all (in this case two) Householder transformations, is the first unit vector.
Since the computed matrix $R$ is upper triangular, this means the
first column of the computed product $QR$ is $[1,~0,~0]^{T}$. Thus, the norm
of the first column of $A-QR$ is exactly $\delta$. There is no reason for $\|A-QR\|$ to
be more than $\delta$, so
the result is that the error $\|A-QR\|$ decreases linearly as $\delta$ drops below $\epsmch^{1/2}$;
although cancellation occurs, the result is to give an increasingly accurate answer
as $\delta$ is reduced. An interesting consequence is that the cancellation apparently cannot result
in arbitrarily poor results; the example illustrated here suggests that, for $A$ with norm one,
$\|A-QR\|$ will perhaps never be significantly greater than $\epsmch^{1/2}$ when 
the wrong sign is used, 
compared to $\epsmch$ when the
correct sign is used (a standard result in numerical linear algebra, 
e.g.\cite[Theorem 19.4]{Hig02}, \cite[Theorem 16.1]{TreBau97}).

\section{History}

According to both Higham \cite{Hig02} and Stewart \cite{Ste98},
the first known use of Householder transformations was
by Turnbull and Aitken in 1932.
Stewart writes ``Householder,
who discovered the transformations independently [in 1958], was the first to realize
their computational significance.'' Stewart also writes ``Householder seems to have missed 
the fact that there are two transformations that
will reduce a vector to a multiple of [the first unit vector] and that the natural construction of one of
them is unstable. This oversight was corrected by Wilkinson [in 1960].'' In Householder's 
1964 book \cite{Hou64}
he writes ``a singularity would arise with one
choice of sign'' (when the two terms cancel exactly) and hence he recommends the other choice of
sign, but, rather surprisingly, he does not mention possible cancellation. Virtually all later books on
numerical linear algebra focus on the latter
issue, motivating the choice \eqref{right}, but we are not aware of any discussion of the 
``inverted V'' phenomenon discussed here. Nor is there any hint that the error $\|A-QR\|$
may be bounded by approximately $\epsmch^{1/2}$ when $A$ has norm one and the wrong
sign is used. Of course, we are not arguing
that using the wrong sign is acceptable. There is no reason to do so, and indeed, even if 
the worst case error is bounded by $\epsmch^{1/2}$, this is still unacceptable when
using the correct sign results in a perfectly stable algorithm.
\bibliography{refs}
\bibliographystyle{alpha}
\end{document}